\documentclass[12pt]{article}
\usepackage{amsfonts}
\newtheorem{lemma}{Lemma}
\newtheorem{theorem}{Theorem}
\newtheorem{proposition}{Proposition}

\newenvironment{proof}{{\bf Proof}.}{\hfill $\Box$}

\begin{document}


\title{\sc The quadratic Fock functor }

\author{\\Luigi Accardi \& Ameur Dhahri\\\vspace{-2mm}\scriptsize Volterra Center, University of Roma Tor Vergata\\\vspace{-2mm}\scriptsize Via Columbia 2, 00133 Roma, Italy
\\\vspace{-2mm}\scriptsize e-mail:accardi@volterra.uniroma2.it\\\scriptsize ameur@volterra.uniroma2.it}
\date{}
\maketitle
\begin{abstract}
We construct the quadratic analogue of the boson Fock functor.
While in the first order (linear) case all contractions on the $1$--particle space
can be second quantized, the semigroup of contractions that admit a 
quadratic second quantization is much smaller due to the nonlinearity.
The encouraging fact is that it contains, as proper sub-groups (i.e. 
the contractions), all the gauge transformations of second kind and all 
the a.e. invertible maps of $\mathbb R^d$ into itself leaving the Lebesgue 
measure quasi-invariant (in particular {\bf all diffeomorphism of $\mathbb R^d$}). 
This allows quadratic $2$-d quantization of gauge theories, of representations of 
the Witt group (in fact it continuous analogue), of the Zamolodchikov hierarchy, 
and much more\dots .
Within this semigroup we characterize the unitary and the isometric elements 
and we single out a class of natural contractions.
\end{abstract}
\section{Introduction}\label{Intr}

The boson (this specification will be omitted in the following) 
Fock functor has its origins in Heisenberg commutation relations.
If $H$ is a complex Hilbert space the Heisenberg $*$--Lie algebra 
$Heis(H)$ is defined by generators.
$$
\{A_g,A_f^+,1 \ \hbox{(central element)} \ : \ f \in H\}
$$
commutation relations
$$
[A_f,A^+_g]=\langle f,g\rangle\cdot 1  \qquad ;\qquad  f,g \in H
$$
(the omitted commutation relations are zero) and involution
$$
(A_f)^\ast=A^+_f\qquad ;\qquad  f\in H
$$
On the universal enveloping algebra of $Heis(H)$, denoted $U(Heis(H))$, 
there is a unique state satisfying
$$
\varphi(1)=1
$$
$$
\varphi(xA_g)=0 \qquad ;\qquad  \forall x \in U(Heis(H)) \ ; \ \forall g\in H
$$
Denoting $\Gamma(H)$ the $GNS$ space of $U(Heis(H))$ with respect to $\varphi$,
the map\\
 $H\mapsto\Gamma (H)$ is a functor defined on the category of Hilbert spaces, 
with morphisms given by contractions to the category of infinite dimensional Hilbert spaces with the same morphisms.\\
$\Gamma(H)$ is called the Fock space over $H$ and, if $V$ is a contraction on $H$ its image $\Gamma (V)$ is called the Fock second quantization of $V$.\\ 
The domain of $\Gamma$ is maximal in the sense that, if $V$ is not a contraction 
on $H$, then $\Gamma(V)$ cannot be a bounded operator on $\Gamma(H)$.

Our goal in this paper is to extend the picture described above, from the Heisenberg algebra, describing the white noise commutation relations, to the algebra describing the commutation relations of the renormalized square of white noise.\\
The algebra of the renormalized square of white noise (RSWN) with test function algebra
$$
{\cal A} := L^2(\Bbb R^d)\cap L^\infty(\Bbb R^d)
$$
is the $*$-Lie-algebra, with central element denoted $1$, generators 
$$
\{ B^+_f,B_h, N_g \ : \ f,g,h\in L^2(\Bbb R^d)\cap L^\infty(\Bbb R^d)\}
$$ 
involution 
$$
(B^+_f)^*=B_f \qquad , \qquad N_f^*=N_{\bar{f}}
$$ 
and commutation relations
$$[B_f,B^+_g]=2c\langle f,g\rangle+4N_{\bar fg},\,\;[N_a,B^+_f]=2B^+_{af},\;c>0$$
$$[B^+_f,B^+_g]=[B_f,B_g]=[N_a,N_{a'}]=0$$
for all $a$, $a'$, $f$, $g\in L^2(\Bbb R^d)\cap L^\infty(\Bbb R^d)$
(the theory can be developed for more general Hilbert algebras, but we will deal only
with this case). 
This is a current algebra over $sl(2,\Bbb R)$ with test function algebra ${\cal A} $.
One can prove that, on the universal enveloping algebra $U(RSWN)$ of the $RSWN$ 
algebra, there exists a unique state $\varphi_F$ such that 
$$
\varphi_F(1)=1
$$
$$
\varphi_F(xB_g)=\varphi_F(xN_f)=0\qquad ;\qquad  \forall f,g\in {\cal A} 
\ ; \ \forall x\in U(RSWN)
$$
By analogy with the Heisenberg algebra, it is natural to call this state 
{\it the quadratic Fock state} and the associated $GNS$ space, denoted 
$\Gamma_2({\cal A})$, {\it the quadratic Fock space}.
The Fock representation of the RSWN is characterized by a cyclic vector $\Phi$,
also called {\it vacuum} as in the first order case, satisfying 
$$
B_f\Phi=N_g\Phi=0
$$ 
for all $f,g\in L^2(\Bbb R^d)\cap L^\infty(\Bbb R^d)$. \\
We refer the interested reader to \cite{AFS}, \cite{AAF} for more details.\\
The extensions, to the quadratic case, of the second quantization procedure 
for linear operators on ${\cal A}$ requires the solution of the following two problems:\\
(1) {\it when does a linear operator on ${\cal A}$ induce a linear operator on} 
$\Gamma_2({\cal A})$?\\
(2) In the cases in which the answer to problem (1) is positive,
 when is the induced operator bounded (a contraction, unitary, isometric, $\dots$)?\\
By inspection on the explicit form of the scalar product of the quadratic Fock space 
(see Lemma \ref{fir} below) one is led to conjecture that two classes of 
linear transformations of ${\cal A}$ should induce contractions on $\Gamma_2({\cal A})$:
\begin{enumerate}
\item[(i)] $*$--endomorphisms of the Hilbert algebra ${\cal A}$
\item[(ii)] generalized gauge transformations of the form
$$
f\mapsto e^{\alpha}f \qquad ; \qquad e^{\alpha}f(x) := e^{\alpha (x)}f(x) 
\ ; \ x\in\Bbb R^d
$$
where $\alpha\in\Bbb R^d\to \Bbb C$ is a complex valued Borel function 
with negative real part (the $-\infty$ value is allowed to include functions
with non full support).
\end{enumerate}
One of our main results is that these are essentially all the linear operators 
on ${\cal A}$ which admit a contractive second quantization on the quadratic Fock space.\\
The scheme of the present paper is the following. In section \ref{quad-Foc-spa}, we recall some properties on the quadratic exponential vectors. Moreover, we prove that the quadratic Fock space is an interacting Fock space with scalar product given explicitly. In section \ref{Quadr-sec-quant}, we characterize those operator on the one--particle Hilbert algebra whose quadratic second quantization is isometric (resp. unitary). In section \ref{Quad-2d-q-contr}, we show with a counter--example that even very simple contractions have a second quantization that is not a contraction and we give a sufficient condition for this to happen. We also introduce the natural candidates for the role of quadratic analogue of the free Hamiltonian evolution and of the Ornstein--Uhlenbeck semigroup.

\section{The quadratic Fock space}\label{quad-Foc-spa}

For $n\in\Bbb N$ the quadratic $n$--particle space is the closed linear span 
of the set 
$$
\big\{B^{+n}_f\Phi \ : \ f\in L^2(\Bbb R^d)\cap L^\infty(\Bbb R^d)\}
$$ 
where by definition $B^{+0}_f\Phi=\Phi$, for all
$f\in L^2(\Bbb R^d)\cap L^\infty(\Bbb R^d)$. The quadratic Fock space 
$\Gamma_2(L^2(\Bbb R^d)\cap L^\infty(\Bbb R^d))$ is the orthogonal sum
of all the quadratic $n$--particle spaces.
The quadratic exponential vector with test function\\
 $f\in L^2(\mathbb{R}^d)\cap L^\infty(\mathbb{R}^d)$, if it exists, is defined by
\begin{equation}\label{def-QEV} 
\Psi(f)=\sum_{n\geq0}\frac{B^{+n}_f\Phi}{n!}
\end{equation}
where by definition
\begin{equation}\label{Psi(0)Phi} 
\Psi(0)= B^{+0}_f\Phi = \Phi
\end{equation}
The following theorem was proved in \cite{AcDh}.
\begin{theorem}\label{prop-QEV} 
The quadratic exponential vector $\Psi(f)$ exists if $\|f\|_\infty<\frac{1}{2}$, and does not exists if  $\|f\|_\infty>\frac{1}{2}$. The set of these vectors is linearly independent and total in 
$\Gamma_2(L^2(\Bbb R^d)\cap L^\infty(\Bbb R^d))$.
Furthermore, the scalar product between two exponential vectors, 
$\Psi(f)$ and $\Psi(g)$, is given by
\begin{equation}\label{Form}
\langle \Psi(f),\Psi(g)\rangle
=e^{-\frac{c}{2}\int_{\mathbb{R}^d}\ln(1-4\bar{f}(s)g(s))ds}
\end{equation}
\end{theorem}
The explicit form of the scalar product between two quadratic $n$--particle vectors 
is due to Barhoumi, Ouerdiane, Riahi \cite{BOR}.
Its proof, which we include for completeness, one needs the following preliminary result
which uses the identity, proved in Proposition 1 of \cite{AcDh}. 
This identity will be frequently used in the following:
\begin{eqnarray}\label{Prop1-AD}
||B^{+m}_f\Phi||^2 &=&c\sum_{k=0}^{m-1}2^{2k+1}\frac{m!(m-1)!}{((m-k-1)!)^2}|\|f^{k+1}\|^2_2\|B_f^{+(m-k-1)}\Phi\|^2\nonumber\\
&=& c\sum_{k=1}^{m-1}2^{2k+1}\frac{m!(m-1)!}{((m-k-1)!)^2}|\|f^{k+1}\|^2_2\|B_f^{+(m-k-1)}\Phi\|^2\nonumber\\
&&+2mc\|f\|^2_2\|B^{+(m-1)}_f\Phi\|^2\nonumber\\
&=& c\sum_{k=0}^{m-2}2^{2k+3}\frac{m!(m-1)!}{(((m-1)-k-1)!)^2}\|f^{k+2}\|^2_2
 \|B_f^{+((m-1)-k-1)}\Phi\|^2  \nonumber\\
&&+  2mc\|f\|^2_2\|B^{+(m-1)}_f\Phi\|^2
\end{eqnarray}
\begin{lemma}\label{A}
For all $f,\,g\in L^2(\Bbb R^d)\cap L^\infty(\Bbb R^d)$ such that 
$\|f\|_\infty<\frac{1}{2}$, $\|g\|_\infty<\frac{1}{2}$, one has
\begin{equation}\label{nth-der-QEV}
\langle B^{+n}_f\Phi,B^{+n}_g\Phi\rangle
={n}!\frac{d^{n}}{dt^{n}}\Big|_{t=0}\langle\Psi(tf),\Psi(g)\rangle
\end{equation}
\end{lemma}
\begin{proof}
Let $f,\,g\in L^2(\Bbb R^d)\cap L^\infty(\Bbb R^d)$ such that $\|f\|_\infty<\frac{1}{2}$, $\|g\|_\infty<\frac{1}{2}$. For all $0\leq t\leq1$, one has
$$
\langle\Psi(tf),\Psi(g)\rangle
=\sum_{m\geq0}\frac{t^m}{(m!)^2}\langle B^{+m}_{f}\Phi, B^{+m}_{g}\Phi\rangle
$$
We now prove that, for $0\leq t\leq1$, the above series can be differentiated (in $t$)
term by term. For all $m\geq n$, one has
\begin{eqnarray*}
\frac{d^n}{dt^n}\Big(\frac{t^m}{(m!)^2}\langle B^{+m}_{f}\Phi, B^{+m}_{g}\Phi\rangle\Big)&=&\frac{m!t^{m-n}}{(m!)^2(m-n)!}\langle B^{+m}_{f}\Phi, B^{+m}_{g}\Phi\rangle\\
&=&\frac{t^{m-n}}{m!(m-n)!}\langle B^{+m}_{f}\Phi, B^{+m}_{g}\Phi\rangle
\end{eqnarray*}
So that, for $0\leq t\leq1$
$$
\Big|\frac{d^n}{dt^n}\Big(\frac{t^m}{(m!)^2}\langle B^{+m}_{f}\Phi, B^{+m}_{g}\Phi\rangle\Big)\Big|\leq U_m
:=\frac{1}{m!(m-n)!}\|B^{+m}_{f}\Phi\| \|B^{+m}_{g}\Phi\|
$$
From the identity (\ref{Prop1-AD}) it follows that
\begin{eqnarray*}
&&c\sum_{k=0}^{m-2}2^{2k+3}\frac{m!(m-1)!}{(((m-1)-k-1)!)^2}\|f^{k+2}\|^2_2
 \|B_f^{+((m-1)-k-1)}\Phi\|^2 \\
&&\leq \Big(4m(m-1)\|f\|^2_\infty\Big)\Big[c\sum_{k=0}^{m-2}2^{2k+1}
 \frac{(m-1)!(m-2)!}{(((m-1)-k-1)!)^2}\|f^{k+1}\|^2_2\\
&&\;\;\;\;\;\;\;\;\;\|B_f^{+((m-1)-k-1)}\Phi\|^2\Big]=\Big(4m(m-1)\|f\|^2_\infty\Big)\|B^{+m}_f\Phi\|^2
\end{eqnarray*}
In conclusion 
$$
||B^{+m}_f\Phi||^2\leq\Big[4m(m-1)\|f\|^2_\infty+2m\|f\|^2\Big]\|B^{+(m-1)}_f\Phi\|^2
$$  
Therefore
\begin{eqnarray*}
\|B^{+m}_{f}\Phi\| \|B^{+m}_{g}\Phi\|&\leq&\sqrt{4m(m-1)\|f\|^2_\infty+2m\|f\|^2_2}\\
&&\sqrt{4m(m-1)\|g\|^2_\infty+2m\|g\|_2^2}\|B^{+(m-1)}_f\Phi\|\|B^{+(m-1)}_g\Phi\|
\end{eqnarray*} 
The definition of $U_m$ then implies that
$$
U_m\leq \frac{\sqrt{4m(m-1)\|f\|^2_\infty+2m\|f\|^2_2}
\sqrt{4m(m-1)\|g\|^2_\infty+2m\|g\|^2_2}}{m(m-n)}\;\,U_{m-1}
$$
If $f$ and $g$ are non-vanishing functions, then
$$
\lim_{m\rightarrow\infty}\frac{U_m}{U_{m-1}}\leq 4\|f\|_\infty\|g\|_\infty<1
$$
because $\|f\|_\infty<\frac{1}{2}$, $\|g\|_\infty<\frac{1}{2}$. 
Hence, the series $\sum_mU_m$ converges. This implies that
$$
\frac{d^{n}}{dt^{n}}\langle\Psi(tf),\Psi(g)\rangle=
\sum_{m\geq n}\frac{t^{m-n}}{m!(m-n)!}\langle B^{+m}_{f}\Phi, B^{+m}_{g}\Phi\rangle
$$
Evaluating the derivative at $t=0$, one obtains (\ref{nth-der-QEV}).
\end{proof}
\begin{lemma}\label{fir}
For all $f,\,g\in L^2(\Bbb R^d)\cap L^\infty(\Bbb R^d)$
the following identity holds 
\begin{eqnarray}\label{s}
\langle B^{+n}_f\Phi,B^{+n}_g\Phi\rangle=\sum_{i_1+2i_2+\dots+ki_k=n}\frac{(n!)^22^{2n-1}c^{i_1+\dots+i_k}}{i_1!\dots
 i_k!2^{i_2}\dots k^{i_k}}
 \langle f,g\rangle^{i_1}\langle f^2,g^2\rangle^{i_2}\dots\langle f^k,g^k\rangle^{i_k}
\end{eqnarray} 
\end{lemma}
\begin{proof} The complex linearity of the map $f\mapsto B^{+}_{f}$ implies that,
for all $\lambda_1,\;\lambda_2\in\mathbb{C}$, 
$$
\langle B^{+n}_{\lambda_1f}\Phi,B^{+n}_{\lambda_2g}\Phi\rangle=
\bar{\lambda}_1^n\lambda_2^n\langle B^{+n}_{f}\Phi,B^{+n}_{g}\Phi\rangle
$$
Therefore it will be sufficient to prove the identity (\ref{s}) for all
$f,\,g\in L^2(\Bbb R^d)\cap L^\infty(\Bbb R^d)$ such that 
$\|f\|_\infty ,\|g\|_\infty<\frac{1}{2}$. In this case one has
\begin{eqnarray}\label{x}
\langle B^{+n}_f\Phi,B^{+n}_g\Phi\rangle&=&n!\frac{d^n}{dt^n}\Big|_{t=0}\langle\Psi(tf),\Psi(g)\rangle\nonumber\\
&=&n!\frac{d^n}{dt^n}\Big|_{t=0}\Big(\exp\big(-\langle\log(1-4t\bar{f}g)\rangle\big)\Big)
\end{eqnarray}
where 
$$
\langle\log(1-4t\bar{f}g)\rangle
:=\frac{c}{2}\int_{\mathbb{R}^d}\log\big(1-4t\bar{f}(s)g(s)\big)ds
$$
Denoting $h(t,s):=\log\big(1-4t\bar{f}(s)g(s)\big)$, its $k$--th derivative (in $t$) is
$$
h^{(k)}(t,s)=2^{2k}(k-1)!(\bar{f}(s))^k(g(s))^k(1-4t\bar{f}(s)g(s))^{-k}
$$
Hence, uniformly for $t\leq1$
\begin{eqnarray}\label{y}
|h^{(k)}(t,s)|\leq \frac{2^{2k}(k-1)!|f(s)|^k|g(s)|^k}{(1-4\|f\|_\infty\|g\|_\infty)^k}
\end{eqnarray}
Thus, the left hand side of (\ref{y}) is integrable in $s$ and
$$
\langle h^{(k)}(t)\rangle
=2^{2k}(k-1)!\int_{\mathbb{R}^d}\frac{(\bar{f}(s))^k(g(s))^k}{(1-4t\bar{f}(s)g(s))^k}ds
$$
Putting $t=0$ one finds
\begin{eqnarray}\label{z}
\langle h^{(k)}(0)\rangle=2^{2k}(k-1)!\langle f^k,g^k\rangle
\end{eqnarray}
Combining the identity (cf. Refs \cite{BOR}, \cite{B})  
\begin{eqnarray}\label{w}
\frac{d^n}{dt^n}e^{\varphi(t)}
=\sum_{i_1+2i_2+\dots+ki_k=n}\frac{2^{2n}n!}{i_1!\dots 
i_k!}\Big(\frac{\varphi^{(1)}(t)}{1!}\Big)^{i_1}\dots
\Big(\frac{\varphi^{(k)}(t)}{k!}\Big)^{i_k}e^{\varphi(t)}
\end{eqnarray}
with (\ref{x}), (\ref{z}) and (\ref{w}) one obtains
\begin{eqnarray*}
\langle B^{+n}_f\Phi,B^{+n}_g\Phi\rangle&=&n!\frac{d^n}{dt^n}\Big|_{t=0}\langle\Psi(tf),\Psi(g)\rangle\\
&=&\sum_{i_1+2i_2+\dots+ki_k=n}\frac{n!2^{2n-1}n!c^{i_1+\dots+i_k}}{i_1!\dots
i_k!2^{i_2}\dots k^{i_k}}\langle f,g\rangle^{i_1}\langle f^2,g^2\rangle^{i_2}
\dots\langle f^k,g^k\rangle^{i_k}
\end{eqnarray*}
from which (\ref{s}) follows.
\end{proof}

The following theorem is an immediate consequence of Lemma \ref{fir}.
\begin{theorem} \label{quad-IFS}
There is a natural ismorphism between the quadratic Fock space
$\Gamma_2(L^2(\mathbb{R}^d)\cap L^\infty(\mathbb{R}^d))$ and the interacting Fock space
$\oplus_{n=0}^{\infty}\otimes^{n}_{symm}\{L^2(\mathbb{R}^d) , \langle \cdot,\cdot\rangle_n\}$,
with scalar products:
$$
\langle f^{\otimes n} , g^{\otimes n}\rangle_n
=\sum_{i_1+2i_2+\dots+ki_k=n}\frac{2^{2n-1}(n!)^2c^{i_1+\dots+i_k}}{i_1!\dots
i_k!2^{i_2}\dots k^{i_k}}\langle f,g\rangle^{i_1}\langle f^2,g^2\rangle^{i_2}
\dots\langle f^k,g^k\rangle^{i_k}
$$ 
\end{theorem}

\section{Quadratic second quantization of contractions } 
\label{Quadr-sec-quant}

Let $T$ be a linear operator on $L^2(\Bbb R^d)\cap L^\infty(\Bbb R^d)$.
If the map
\begin{equation}\label{df-G2T}
\Psi(f)\mapsto \Psi(Tf)
\end{equation}
is well defined for all quadratic exponential vectors then, by the linear 
independence of these vectors, it admits a linear extension to a dense subspace 
of $\Gamma_2(L^2(\Bbb R^d)\cap L^\infty(\Bbb R^d))$, denoted $\Gamma_2(T)$
and called {\it the quadratic second quantization of $T$}.\\
From (\ref{Psi(0)Phi}) and (\ref{df-G2T}) it follows 
that, if $\Gamma_2(T)$ exists then, whatever $T$ is, it leaves the quadratic vacuum invariant:
$$
\Gamma_2(T)\Phi = \Phi
$$
\begin{lemma}
Let $T$ be a linear operator on $L^2(\Bbb R^d)\cap L^\infty(\Bbb R^d)$. 
Then $\Gamma_2(T)$ is well defined on the set of all the quadratic exponential vectors 
$$\{\Psi(f)\quad,\qquad f\in  L^2(\Bbb R^d)\cap L^\infty(\Bbb R^d) \mbox{ s.t }\|f\|_\infty<\frac{1}{2}\}$$
if and only if $T$ is a contraction on $L^2(\mathbb{R}^d)\cap L^\infty(\mathbb{R}^d)$
equipped with the norm $\|.\|_\infty$.    
\end{lemma}
\begin{proof}
Sufficiency. If $T:L^\infty(\mathbb{R}^d)\rightarrow L^\infty(\mathbb{R}^d)$ 
is a contraction, then\\ 
$\|Tf\|_\infty \leq \|f\|_\infty < 1/2$ for any test function 
$f\in L^2(\mathbb{R}^d)\cap L^\infty(\mathbb{R}^d)$ such that $\|f\|_\infty < 1/2$. 
Therefore $\Gamma_2(T)\Psi(f)$ is well defined.\\ 
Necessity. If $\Gamma_2(T)$ is well defined, then one has $\|Tg\|_\infty\leq\frac{1}{2}$, 
for any\\ 
$g\in L^2(\Bbb R^d)\cap L^\infty(\Bbb R^d)$ such that $\|g\|_\infty<\frac{1}{2}$.
By linearity $T$ maps the open unit $\|.\|_\infty$--ball of 
$L^2(\Bbb R^d)\cap L^\infty(\Bbb R^d)$ into the closed unit $\|.\|_\infty$--ball of 
$L^2(\Bbb R^d)\cap L^\infty(\Bbb R^d)$ , i.e. it is a contraction.
\end{proof}

\subsection{Isometric and unitarity characterization of the quadratic second quantization}\label{Iso-un-q2d-quan}

Let us start by giving a sufficient condition on $T$, which ensures that 
$\Gamma_2(T)$ is an isometry (resp. unitary operator).

A Hilbert algebra endomorphism (resp. automorphism) $T$ of\\
 $L^2(\mathbb{R}^d)\cap L^\infty (\mathbb{R}^d)$ is said to be 
 a {\it $*$-endomorphism} (resp. {\it $*$-automorphism}) if $T$ is an isometry
 (resp. a unitary operator) with respect to the pre-Hilbert structure of $L^2(\mathbb{R}^d)\cap L^\infty(\mathbb{R}^d)$, which satisfies
$$T(fg)=T(f)T(g),\;\,(T(f))^*=T(\bar{f}).$$ 
The following proposition is an immediate consequence of Lemma \ref{fir}.
\begin{proposition}\label{pro}
If $\alpha:\mathbb{R}^d\rightarrow\mathbb{R}$ is a Borel function, 
$T_1$ is a $*$-endomorphism of $L^2(\Bbb R^d)\cap L^\infty(\Bbb R^d)$ and 
$$
T :=e^{i\alpha}T_1
$$ 
then $\Gamma_2(T)$ is an isometry. Moreover, if $T_1$ is a $*$-automorphism of\\ 
$L^2(\mathbb{R}^d)\cap L^\infty(\mathbb{R}^d)$, then $\Gamma_2(T)$ is unitary.
\end{proposition}
\begin{proof}
To prove that $\Gamma_2(T)$ is an isometry it is sufficient to prove that
it preserves the scalar product of two arbitray quadratic exponential vectors.
From (\ref{def-QEV}) and the mutual orthogonality of different $n$--particle spaces,
it will be sufficient to prove that, for each $n\in\mathbb N$ and
$f,g\in L^2(\Bbb R^d)\cap L^\infty(\Bbb R^d)$ one has:
$$
\langle B^{+n}_{Tf}\Phi,B^{+n}_{Tg}\Phi\rangle 
 = \langle B^{+n}_{f}\Phi,B^{+n}_{g}\Phi\rangle
$$
and, because of Lemma \ref{fir}, this identity follows from
$$
\langle (Tf)^k,(Tg)^k\rangle = \langle f^k,g^k\rangle
\qquad ; \qquad \forall k\in\mathbb N
\qquad ; \qquad \forall f,g\in L^2(\Bbb R^d)\cap L^\infty(\Bbb R^d)
$$
But this identity holds because our assumptions on $T$ imply that
$$
\langle (Tf)^k,(Tg)^k\rangle 
= \langle e^{ik\alpha}(T_1f)^k,e^{ik\alpha}(T_1g)^k\rangle
= \langle T_1(f^k),T_1(g^k)\rangle
= \langle f^k,g^k\rangle
$$
Thus $\Gamma_2(T)$ is an isometry. If, in addition, $T_1$ is a $*$-automorphism of\\ 
$L^2(\mathbb{R}^d)\cap L^\infty(\mathbb{R}^d)$, then $T$ is surjective. Hence
the range of $\Gamma_2(T)$, containing all the quadratic exponential vectors, 
is the whole quadratic Fock space. The thesis then follows because an isometry
with full range is unitary.
\end{proof}

In the following our goal is to prove the converse of the above proposition.
\begin{lemma}\label{B}
\begin{enumerate}
\item[i)] If $\Gamma_2(T)$ is a unitary operator, then 
\begin{eqnarray}\label{tun}
\langle(Tf)^{n},(Tg)^{n}\rangle=\langle f^{n},g^{n}\rangle
\end{eqnarray}
for all $n\in\mathbb{N}^*$ and $f,\,g\in  L^2(\Bbb R^d)\cap L^\infty(\Bbb R^d)$. 
\item[ii)] If $\Gamma_2(T)$ is an isometry, then for all $n\in\mathbb{N}^*$ and $f\in  L^2(\Bbb R^d)\cap L^\infty(\Bbb R^d)$ 
$$\|(Tf)^{n}\|_2=\|f^{n}\|_2$$
\end{enumerate}
\end{lemma}
\begin{proof}
Suppose that $\Gamma_2(T)$ is a unitary operator. Let us fix two functions $f,\,g\in  L^2(\Bbb R^d)\cap L^\infty(\Bbb R^d)$ such that $\|f\|_\infty<\frac{1}{2},\,\|g\|_ \infty<\frac{1}{2}$. Then, one has
$$\langle\Psi(Tf),\Psi(Tg)\rangle=\langle \Psi(f),\Psi(g)\rangle$$
 It follows that 
$$\langle\Psi(tTf),\Psi(Tg)\rangle=\langle \Psi(tf),\Psi(g)\rangle$$
for all $t$ such that $|t|<1$. Therefore, Lemma \ref{A} implies that
\begin{eqnarray}\label{samiha}
\langle B^{+n}_{Tf}\Phi,B^{+n}_{Tg}\Phi\rangle=\langle B^{+n}_{f}\Phi,B^{+n}_{g}\Phi\rangle
\end{eqnarray}
for all $n\in\mathbb{N}$. Let us prove the statement i) by induction. 

- For $n=1$, we have
$$\langle B^{+}_{Tf}\Phi,B^{+}_{Tg}\Phi\rangle=\langle B^{+}_{f}\Phi,B^{+}_{g}\Phi\rangle$$
This gives
$$
\langle Tf,Tg\rangle=\langle f,g\rangle
$$
- Suppose that (\ref{tun}) holds for $k\leq n$. 
Then, from (\ref{samiha}) and  the identity (\ref{Prop1-AD}), one obtains
\begin{eqnarray*}
&&\langle B^{+(n+1)}_{Tf}\Phi,B^{+(n+1)}_{Tg}\Phi\rangle\\
&&\;\;=c\sum^{n}_{k=0}2^{2k+1}{n!(n+1)!\over((n-k)!)^2}\langle (Tf)^{k+1}, (Tg)^{k+1}\rangle\langle B^{+(n-k)}_{Tf}\Phi,B^{+(n-k)}_{Tg}\Phi\rangle\\
&&\;\;=2^{2n+1}n!(n+1)!c\,\langle (Tf)^{n+1}, (Tg)^{n+1}\rangle\\
&&\;\;\;\;+c\sum^{n-1}_{k=0}2^{2k+1}{n!(n+1)!\over((n-k)!)^2}\langle (Tf)^{k+1}, (Tg)^{k+1}\rangle\langle B^{+(n-k)}_{Tf}\Phi,B^{+(n-k)}_{Tg}\Phi\rangle\\
&&\;\;=2^{2n+1}n!(n+1)!c\,\langle f^{n+1}, g^{n+1}\rangle\\
&&\;\;\;\;+c\sum^{n-1}_{k=0}2^{2k+1}{n!(n+1)!\over((n-k)!)^2}\,\langle f^{k+1}, g^{k+1}\rangle\langle B^{+(n-k)}_{f}\Phi,B^{+(n-k)}_g\Phi\rangle
\end{eqnarray*}
By the induction assumption, one has
\begin{eqnarray*}
&&c\sum^{n-1}_{k=0}2^{2k+1}{n!(n+1)!\over((n-k)!)^2}\,
  \langle (Tf)^{k+1}, (Tg)^{k+1}\rangle\langle B^{+(n-k)}_{Tf}\Phi,B^{+(n-k)}_{Tg}\Phi\rangle \\
&&\;\;=c\sum^{n-1}_{k=0}2^{2k+1}{n!(n+1)!\over((n-k)!)^2}\,\langle f^{k+1}, g^{k+1}\rangle\langle B^{+(n-k)}_{f}\Phi,B^{+(n-k)}_g\Phi\rangle
\end{eqnarray*}
which implies that
$$
\langle(Tf)^{n+1},(Tg)^{n+1}\rangle=\langle f^{n+1},g^{n+1}\rangle
\qquad ; \qquad \forall n\in\mathbb{N}^*
$$
Thus (\ref{tun}) holds for all $n\in\mathbb{N}^*$.\\
The proof of statement ii) is obtained by replacing, in the above argument,
the test function $g$ by $f$.
\end{proof}
\begin{lemma}\label{C}
Suppose that $\Gamma_2(T)$ is an isometry. Then, for any $I\subset\mathbb{R}^d$
such that $|I|<\infty$, one has
$$
|T(\chi_I)(x)|=1
$$
on $supp(T(\chi_I))$ a.e.
\end{lemma}
\begin{proof}
By assumption $\Gamma_2(T)$ is an isometry, hence from Lemma \ref{B},
$\forall n\in\mathbb N$:
\begin{eqnarray}\label{iso-pres-meas}
\langle (T(\chi_I))^n, (T(\chi_I))^n\rangle=\langle (\chi_I)^n, (\chi_I)^n\rangle
=\langle \chi_I, \chi_I\rangle=|I|
\end{eqnarray}
for any subset $I\subset\mathbb{R}^d$ such that $|I|<\infty$.  But, one has
\begin{eqnarray}\label{amal}
\langle (T(\chi_I))^n, (T(\chi_I))^n\rangle=
\big|\big\{x\in\mathbb{R}^d,\,|T(\chi_I)(x)|=1\big\}\big|+\int_J|T(\chi_I)(x)|^{2n}dx
\end{eqnarray}
where $\big| \ \cdot \ \big|$ denotes Lebesgue measure and
$$
J:=\{x\in\mathbb{R}^d,\,|T(\chi_I)(x)|\neq1\;\mbox{ and }\;|T(\chi_I)(x)|>0\}
$$
Since the identity (\ref{amal}) holds $\forall n\in\mathbb N$, it follows that
$$
\int_J|T(\chi_I(x)|^{2n}dx=\int_J|T(\chi_I(x)|^{2(n+1)}dx
\qquad ; \qquad \forall n\in\mathbb N
$$
But it is not difficult to prove that this is impossible if $|J|>0$.
\end{proof}
\begin{lemma}\label{D}
If $I\subset\mathbb{R}^d$ such that $|I|<\infty$ and $\Gamma_2(T)$ is an isometry, 
then there exist a function $\alpha_I:\mathbb{R}^d\rightarrow\mathbb{R}$ 
and a subset $\tau(I)\subset\mathbb{R}^d$ such that 
$$
T(\chi_I)=e^{i\alpha_I}\chi_{\tau(I)}
$$ 
and $|I|=|\tau(I)|$. Moreover,  if $I_1,\,I_2$ is an arbitrary partition  of $I$, then
\begin{eqnarray}\label{add-tau}
\tau(I)=\tau(I_1)\cup\tau(I_2) \qquad , \qquad  a.e. 
\end{eqnarray}
In particular, if $I_1\subset I$, then a.e. $\tau(I_1)\subset\tau(I)$.
\end{lemma}
\begin{proof}
Lemma \ref{C} implies that there exist a function $\alpha_I:\mathbb{R}^d\rightarrow\mathbb{R}$ and a subset $\tau(I)\subset\mathbb{R}^d$ such that $T(\chi_I)=e^{i\alpha_I}\chi_{\tau(I)}$. From (\ref{iso-pres-meas}) one has
$$
|\tau(I)|=\langle T(\chi_I),T(\chi_I)\rangle=\langle \chi_I,\chi_I\rangle=|I|
$$
Let $I_1,\,I_2$ be a partition of $I$. From 
$\chi_{I}=\chi_{I_1\cup I_2}=\chi_{I_1}+\chi_{I_2}$, it follows that
$$
T(\chi_I)= T(\chi_{I_1})+ T(\chi_{I_2})
$$
i.e.
$$
e^{i\alpha_I}\chi_{\tau(I)}
=e^{i\alpha_{I_1}}\chi_{\tau(I_1)}+e^{i\alpha_{I_2}}\chi_{\tau(I_2)}
$$
Multiplying both sides by $\chi_{\tau(I_1)\cup \tau(I_2)}$, one finds
$$
e^{i\alpha_I}\chi_{\tau(I)\cap[\tau(I_1)\cup\tau(I_2)]} = e^{i\alpha_{I_1}}\chi_{\tau(I_1)}+e^{i\alpha_{I_2}}\chi_{\tau(I_2)}\\
 =  e^{i\alpha_I}\chi_{\tau(I)}
$$    
Therefore, one has $\tau(I)=\tau(I_1)\cup\tau(I_2)$ a.e. 
Since the partition $I_1,\,I_2$ of $I$ is arbitrary, it follows that 
$I_1\subset I$ implies that $\tau(I_1)\subset\tau(I).$
\end{proof}
\begin{lemma}\label{E}
If $\Gamma_2(T)$ is an isometry and $I_1,\,I_2\subset\mathbb{R}^d$ are such that\\
 $|I_1|<\infty,\,|I_2|<\infty$ and $|I_1\cap I_2|=0$, then $|\tau(I_1)\cap \tau(I_2)|=0.$
\end{lemma}
\begin{proof}
Suppose that $|I_1\cap I_2|=0$. Then, from the identity 
$$\chi_{I_1\cup I_2}=\chi_{I_1}+\chi_{I_2}-\chi_{I_1\cap I_2}$$
it follows that, a.e.
$$
\chi_{I_1\cup I_2}=\chi_{I_1}+\chi_{I_2}
$$
and therefore also 
$$
T(\chi_{I_1\cup I_2})=T(\chi_{I_1})+T(\chi_{I_2})\qquad ; \qquad  a.e
$$
Applying (\ref{iso-pres-meas}) one then gets
\begin{eqnarray}\label{ism}
|I_1|+|I_2|&=&\langle \chi_{I_1\cup I_2},\chi_{I_1\cup I_2}\rangle\nonumber\\
&=&\langle T(\chi_{I_1\cup I_2}),T(\chi_{I_1\cup I_2})\rangle\nonumber\\
&=&\langle T(\chi_{I_1}),T(\chi_{I_1})\rangle+\langle T(\chi_{I_2}),T(\chi_{I_2})\rangle\nonumber\\
&&+\langle T(\chi_{I_1}),T(\chi_{I_2})\rangle+\langle T(\chi_{I_2}),T(\chi_{I_1})\rangle\nonumber\\
&=&|I_1|+|I_2|+\int_{\tau(I_1)\cap\tau(I_2)}e^{i(\alpha_{I_2}-\alpha_{I_1})}(x)dx\nonumber\\
&&+\int_{\tau(I_1)\cap\tau(I_2)}e^{-i(\alpha_{I_2}-\alpha_{I_1})}(x)dx\nonumber\\
&=&|I_1|+|I_2|+2\int_{\tau(I_1)\cap\tau(I_2)} cos((\alpha_{I_2}-\alpha_{I_1})(x))dx
\end{eqnarray}    
which implies that
\begin{eqnarray}\label{nc}
\int_{\tau(I_1)\cap\tau(I_2)} cos((\alpha_{I_2}-\alpha_{I_1})(x))dx=0
\end{eqnarray}
Put $I=I_1\cup I_2$. From the identities
\begin{eqnarray*}
e^{i\alpha_I}\chi_{\tau(I)}&=&e^{i\alpha_{I_1}}\chi_{\tau(I_1)}+e^{i\alpha_{I_2}}\chi_{\tau(I_2)}\\
\tau(I)&=&\tau(I_1)\cup\tau(I_2)\qquad  a.e
\end{eqnarray*}
it follows that if $x\in \tau(I_1)\cap\tau(I_2)$, then
$$e^{i\alpha_I}(x)=e^{i\alpha_{I_1}}(x)+e^{i\alpha_{I_2}}(x)$$
Thus, one obtains
$$e^{i(\alpha_{I}(x)-\alpha_{I_1}(x)}=1+e^{i(\alpha_{I_2}(x)-\alpha_{I_1}(x))}$$
This gives
$$1=|1+e^{i(\alpha_{I_2}(x)-\alpha_{I_1}(x))}|^2=2+2cos(\alpha_{I_2}(x)-\alpha_{I_1}(x))$$
which yields that
$$cos(\alpha_{I_2}(x)-\alpha_{I_1}(x))=-\frac{1}{2}$$
This, together with (\ref{nc}) implies that $|\tau(I_1)\cap\tau(I_2)|=0.$
\end{proof}
\begin{lemma}\label{F}
In the notations and assumptions of Lemma \ref{D}, for any $I\subset\mathbb{R}^d$ 
such that $|I|<\infty$ and any $I_1\subset I$ one has
$$
e^{i\alpha_{I_1}}\chi_{\tau(I_1)}=e^{i\alpha_I}\chi_{\tau(I_1)}
$$
for almost any $x\in \tau(I_1).$ 
\end{lemma}
\begin{proof}
Let $I_2=I\setminus I_1$.  Arguing as in the proof of of Lemma \ref{D} one finds that 
$$
e^{i\alpha_I}\chi_{\tau(I)}
=e^{i\alpha_{I_1}}\chi_{\tau(I_1)}+e^{i\alpha_{I_2}}\chi_{\tau(I_2)}
$$
Thus, if we multiply the two sides in the above identity by $\chi_{\tau(I_1)}$, then from Lemmas \ref{D}, \ref{E}, it follows that
$$
e^{i\alpha_I}\chi_{\tau(I_1)}=e^{i\alpha_{I_1}}\chi_{\tau(I_1)} \qquad , \qquad a.e
$$  
\end{proof}

\begin{lemma}\label{G}
In the notations and assumptions of Lemma \ref{D} there exists a function
$\alpha:\mathbb{R}^d\rightarrow\mathbb{R}$ such that for any $I\subset\mathbb{R}^d$,
with $|I|<\infty$ 
$$
T(\chi_I)=e^{i\alpha}\chi_{\tau(I)}
$$
where $\tau(I)\subset \mathbb{R}^d$ and $|\tau(I)|=|I|.$
\end{lemma}
\begin{proof}
Let $(I_n)_n$ be an increasing sequence of subsets of $\mathbb{R}^d$ such that\\ $|I_n|<\infty,\,\forall n\in \mathbb{N}$ and $\bigcup_{n\in\mathbb{N}}I_n=\mathbb{R}^d$.
Define the function $\alpha:\mathbb{R}^d \to \mathbb{R}$ by $\alpha(x)=\alpha_{I_n}(x),$ 
for any $n\in\mathbb{N}$ such that $x\in I_n$, where $\alpha_{I_n}$ is defined
as in Lemma (\ref{D}). Then $\alpha$ is well defined because, denoting
$$
n(x):=\min\{n\in\mathbb{N},\,x\in I_n\}
\qquad ; \qquad x\in\mathbb{R}^d
$$
Lemma \ref{F} implies that, for any $m,\,n\in\mathbb{N}$ such that $n(x)\leq m\leq n$, 
$$
e^{i\alpha_{I_m}}\chi_{\tau(I_m)}=e^{i\alpha_{I_n}}\chi_{\tau(I_m)}
$$
In particular, for any $n\geq n(x)$, one has
$$e^{i\alpha_{I_n}}\chi_{\tau(I_{n(x)})}=e^{i\alpha_{I_{n(x)}}}\chi_{\tau(I_{n(x)})}$$
which implies that
$$
\alpha_{I_n}(x)=\alpha_{I_n(x)},\;\forall n\geq n(x)
$$
This ends the proof of the above lemma.
\end{proof}

Using all together Proposition \ref{pro}, Lemmas \ref{B}, \ref{D} and \ref{G}, we prove the following.
\begin{theorem}\label{struc-iso-un}
$\Gamma_2(T)$ is an isometry (resp. unitary) if and only if there exist a function $\alpha$ from $\mathbb{R}^d$ to $\mathbb{R}$ and a $*$-endomorphism (resp. $*$-automorphism) $T_1$ of $L^2(\mathbb{R}^d)\cap L^\infty(\mathbb{R}^d)$ such that
$$T=e^{i\alpha}T_1$$
\end{theorem}
\begin{proof}
Sufficiency has been proved in Proposition \ref{pro}. \\
Necessity. Suppose that $\Gamma_2(T)$ is an isometry. Then, from Lemma \ref{B}, 
$T$ is an isometry. Moreover, Lemma \ref{G} implies that there exists a function
$\alpha:\mathbb{R}^d\rightarrow\mathbb{R}$ such that for any $I\subset\mathbb{R}^d$,
$|I|<\infty$
$$
T(\chi_I)=e^{i\alpha}\chi_{\tau(I)}
$$
where $\tau(I)\subset \mathbb{R}^d$ and $|\tau(I)|=|I|$. Define 
the map $T_1$ by:
\begin{eqnarray}\label{df-T1}
T_1:\chi_I\in L^2(\mathbb{R}^d)\cap L^\infty(\mathbb{R}^d) \to
T_1(\chi_I):=\chi_{\tau(I)}
\end{eqnarray}
for all $I\subset \mathbb{R}^d$ such that $|I|<\infty$. 
In order to prove that $T_1$ extends, by linearity and continuity, to a 
$*$-endomorphism of $L^2(\mathbb{R}^d)\cap L^\infty(\mathbb{R}^d)$, 
it is sufficient to prove that for all $I,J\subset \mathbb R$ with $|I|<\infty,\,|J|<\infty$
\begin{eqnarray}\label{Cha}
T_1(\chi_I\chi_J)=T_1(\chi_I)T_1(\chi_J)
=\chi_{\tau(I)}\chi_{\tau(J)} =\chi_{\tau(I)\cap \tau(J)}
\qquad ,\qquad  a.e
\end{eqnarray}
But, by definition of $T_1$ one has
$$
T_1(\chi_I\chi_J)=T_1(\chi_{I\cap J})=\chi_{\tau(I\cap J)}
$$
therefore our thesis is equivalent to
\begin{eqnarray}\label{Cha1}
\tau(I)\cap \tau(I) = \tau(I\cap J)
\qquad ,\qquad  a.e
\end{eqnarray}
Finally, since from Lemma \ref{D} we know that $\tau(I\cap J)\subset \tau(I)\cap\tau(J)$,
(\ref{Cha1}) will follow if we prove that
\begin{eqnarray}\label{Cha2}
|\tau(I)\cap\tau(J)|=|\tau(I\cap J)|
\end{eqnarray}
To prove (\ref{Cha2}) notice that, since $T$, hence $T_1$, is an isometry, one has
\begin{eqnarray}\label{lui}
\langle T_1(\chi_{I\cup J}),T_1(\chi_{I\cup J})\rangle
=\langle\chi_{I\cup J},\chi_{I\cup J}\rangle
=|I|+|J|-|I\cap J|
\end{eqnarray}
On the other hand, from Lemma \ref{D} we know that the map $I\mapsto\tau(I)$
is finitely addditive, hence monotone. Therefore, using linearity, (\ref{df-T1})
and the identity $\chi_{I\cup J}= \chi_{I}+\chi_{J}-\chi_{I\cap J}$, we find
\begin{eqnarray*}
\langle T_1(\chi_{I\cup J}),T_1(\chi_{I\cup J})\rangle&=&\langle T_1(\chi_{I})+T_1(\chi_{J})-T_1(\chi_{I\cap J}),T_1(\chi_{I})\\
&&\;\;\;\;\;\;\;+T_1(\chi_{J})-T_1(\chi_{I\cap J})\rangle\\
&=&\langle T_1(\chi_{I}),T_1(\chi_{I})\rangle+\langle T_1(\chi_{I}),T_1(\chi_{J})\rangle\\
&&-\langle T_1(\chi_{I}),T_1(\chi_{I\cap J})\rangle+\langle T_1(\chi_{J}),T_1(\chi_{I})\rangle\\
&&+\langle T_1(\chi_{J}),T_1(\chi_{J})\rangle-\langle T_1(\chi_{J}),T_1(\chi_{I\cap J})\rangle\\
&&-\langle T_1(\chi_{I\cap J}),T_1(\chi_{I})\rangle-\langle T_1(\chi_{I\cap J}),T_1(\chi_{J})\rangle\\
&&+\langle T_1(\chi_{I\cap J}),T_1(\chi_{I\cap J})\rangle\\
&=&\langle  \chi_{\tau(I)}, \chi_{\tau(I)}\rangle+\langle  \chi_{\tau(I)}, \chi_{\tau(J)}\rangle-\langle  \chi_{\tau(I)}, \chi_{\tau(I\cap J)}\rangle\\
&&+\langle  \chi_{\tau(J)}, \chi_{\tau(I)}\rangle+\langle  \chi_{\tau(J)}, \chi_{\tau(J)}\rangle-\langle  \chi_{\tau(J)}, \chi_{\tau(I\cap J)}\rangle\\
&&-\langle  \chi_{\tau(I\cap J)}, \chi_{\tau(I)}\rangle-\langle  \chi_{\tau(I\cap J)}, \chi_{\tau(J)}\rangle\\
&&+\langle  \chi_{\tau(I\cap J)}, \chi_{\tau(I\cap J)}\rangle
\end{eqnarray*}
Using the isometry property and the fact that $\tau(I\cap J)\subseteq \tau(I)\cap\tau(J)$,
we see that this expression is equal to 
\begin{eqnarray*}
&&|I| + |\tau(I)\cap \tau(J)| -  |\tau(I\cap J)|+ |\tau(I)\cap \tau(J)| + |J|-  |\tau(I\cap J)|\\
&&\;\;\;\;\;\;-  |\tau(I\cap J)|-  |\tau(I\cap J)| +  |\tau(I\cap J)|\\
&&\;\;\;=  |I| + |J| + 2|\tau(I)\cap \tau(J)| -  3|\tau(I\cap J)|\\
&&\;\;\;=  |I| + |J| + 2|\tau(I)\cap \tau(J)| -  3|I\cap J|
\end{eqnarray*}
Since this is equal to the right hand side of (\ref{lui}), 
we conclude that
$$
-|I\cap J| = 2|\tau(I)\cap \tau(J)| -  3|I\cap J|\Leftrightarrow
|\tau(I)\cap \tau(J)|=|I\cap J|=|\tau(I\cap J)|
$$
which is equivalent to (\ref{Cha2}) and therefore to (\ref{Cha}).\\
Since a unitary operator is an isometry we conclude that, if $\Gamma_2(T)$ is unitary,
then $T_1$, defined by (\ref{df-T1}), is an invertible $*$-endomorphism of
$L^2(\mathbb{R}^d)\cap L^\infty(\mathbb{R}^d)$, i.e. a $*$-automorphism.  
\end{proof}

\section{Quadratic second quantization of contractions }
\label{Quad-2d-q-contr}

We will use the following remark.\\
{\bf Remark} 
Let $A=(a_{ij})_{i,j},\;B=(b_{ij})_{i,j},\;C=(c_{ij})_{i,j}$ and 
$D=(d_{ij})_{i,j}$ be matrices such that, in the operator order:
$$
0\leq A\leq B,\;\mbox{ and }\,\,0\leq C\leq D
$$
Then by Schur's Lemma
$$
0\leq \left((b_{ij} - a_{ij})c_{ij}\right)_{i,j}\Leftrightarrow
\left(a_{ij}c_{ij}\right)_{i,j}\leq \left(b_{ij} c_{ij}\right)_{i,j}
$$
$$
0\leq \left(b_{ij}(d_{ij} - c_{ij})\right)_{i,j}\Leftrightarrow
\left(b_{ij}c_{ij}\right)_{i,j}\leq \left(b_{ij} d_{ij}\right)_{i,j}
$$
Consequently one has 
\begin{eqnarray}\label{Kernel}
0\leq (a_{ij}c_{ij})_{i,j}\leq(b_{ij}d_{ij})_{i,j}
\end{eqnarray}

\begin{theorem}\label{th-GSGT=GST}
The set of all operators $T$ such that $\Gamma_2(T)$ is a contraction is a multiplicative semigroups denoted $\hbox{Contr}_2(L^2\cap L^\infty)$. Moreover
\begin{equation}\label{GSGT=GST} 
\Gamma_2(S)\Gamma_2(T)=\Gamma_2(ST)\qquad ; \ \forall S,T\in\hbox{Contr}_2(L^2\cap L^\infty)
\end{equation}
\end{theorem}
\begin{proof} Let $S,T\in\hbox{Contr}_2(L^2\cap L^\infty)$.
Then $\Gamma_2(S)$, $\Gamma_2(T)$ and hence $\Gamma_2(S)\Gamma_2(T)$ is a contraction on $\Gamma_2(L^2\cap L^\infty)$.
Therefore it is uniquely determined by its value on the quadratic exponential vectors. If $\Psi(f)$ is such a vector, then
$$
\Gamma_2(S)\Gamma_2(T)\Psi(f)=
\Gamma_2(S)\Psi(Tf)=\Psi(STf)=\Gamma_2(ST)\Psi(f)
$$
Thus $\Gamma_2(ST)$ is a contraction and (\ref{GSGT=GST}) holds.
\end{proof}

Now, we prove the following.
\begin{proposition}\label{spec-contr}
If $T=\mathcal{M}_\varphi T_1$ is a contraction for $L^2(\mathbb{R}^d)$ and $L^\infty(\mathbb{R}^d)$, where $\varphi\in L^2(\mathbb{R}^d)\cap L^\infty(\mathbb{R}^d)$ such that $\|\varphi\|_\infty\leq 1$ and $T_1$ is an homomorphism of $L^2(\mathbb{R}^d)\cap L^\infty(\mathbb{R}^d)$, then $\Gamma_2(T)$ is a contraction.
\end{proposition}
\begin{proof} We have
\begin{eqnarray}\label{chah}
\|\Gamma_2(T)\big(\alpha_1\Psi(f_1)+\dots+\alpha_l\Psi(f_l)\big)\|^2&=&\|\alpha_1\Psi(Tf_1)+\dots+\alpha_l\Psi(Tf_l)\|^2\nonumber\\
&=&\sum_{i,j=1}^l\bar{\alpha}_i\alpha_j\langle \Psi(Tf_i),\Psi(Tf_j)\rangle\\
&=&\sum_{n\geq0}\frac{1}{(n!)^2}\Big[\sum_{i,j=1}^l\bar{\alpha}_i\alpha_j\langle B^{+n}_{Tf_i}\Phi,B^{+n}_{Tf_j}\Phi\rangle\Big].\nonumber
\end{eqnarray}
Put
$$A_{n,T}=\Big(\langle B^{+n}_{Tf_i}\Phi,B^{+n}_{Tf_j}\Phi\rangle\Big)_{i,j},\;A_{n}=\Big(\langle B^{+n}_{f_i}\Phi,B^{+n}_{f_j}\Phi\rangle\Big)_{i,j}.$$
Now, our purpose is to prove, under the assumptions of the above proposition, that
\begin{eqnarray}\label{Aaa}
0\leq A_{n,T}\leq A_n,
\end{eqnarray}
for all $n\in\mathbb{N}$.

Note that, for $v=(\alpha_1,\dots,\alpha_l)$, one has
\begin{eqnarray*}
\langle v,A_{n,T}v\rangle&=&\sum_{i,j=1}^l\bar{\alpha}_i\alpha_j \langle B^{+n}_{Tf_i}\Phi,B^{+n}_{Tf_j}\Phi\rangle\\
&=&\|\alpha_1 B^{+n}_{Tf_1}\Phi+\dots+\alpha_l B^{+n}_{Tf_l}\Phi\|^2.
\end{eqnarray*}
This implies that $A_{n,T}$ is a positive matrix. Now, let us prove the second inequality in (\ref{Aaa}) by induction on $n$.

- For $n=1$, one has
\begin{eqnarray*}
\langle v,A_{1,T}v\rangle&=&\sum_{i,j=1}^l\bar{\alpha}_i\alpha_j \langle B^{+}_{Tf_i}\Phi,B^{+}_{Tf_j}\Phi\rangle\\
&=&2c \sum_{i,j=1}^l\bar{\alpha}_i\alpha_j \langle Tf_i,Tf_j\rangle\\
&=&2c\|T\big(\alpha_1f_1+\dots+\alpha_lf_l\big)\|^2_2\\
&\leq&2c\|\alpha_1f_1+\dots+\alpha_lf_l\|^2_2.
\end{eqnarray*}
Because
$$2c\|\alpha_1f_1+\dots+\alpha_lf_l\|^2_2=\sum_{i,j=1}^l\bar{\alpha}_i\alpha_j \langle B^{+}_{f_i}\Phi,B^{+}_{f_j}\Phi\rangle=\langle v,A_{1}v\rangle,$$
one obtains that $A_{1,T}\leq A_1$.

- Let $n\geq 1$ and suppose that $A_{n,T}\leq A_n$. Note that for any $f,\,g\in L^2(\mathbb{R}^d)\cap L^\infty(\mathbb{R}^d)$, Proposition 1 of \cite{AcDh} implies that
\begin{eqnarray*}
\langle B^{+(n+1)}_f\Phi,B^{+(n+1)}_g\Phi\rangle&=&c\sum^{n}_{k=0}2^{2k+1}{n!(n+1)!\over((n-k)!)^2}\,\langle f^{k+1}, g^{k+1}\rangle\\
\;\;\;\;\;&&\langle B^{+(n-k)}_f\Phi,B^{+(n-k)}_g\Phi\rangle.
\end{eqnarray*}
Then, one gets
\begin{eqnarray*}
\langle v,A_{n+1,T}v\rangle&=&\sum_{i,j=1}^l\bar{\alpha}_i\alpha_j\langle B^{+(n+1)}_{Tf_i}\Phi,B^{+(n+1)}_{Tf_j}\Phi\rangle\\
&=&c\sum_{k=0}^n2^{2k+1}\frac{(n+1)!n!}{((n-k)!)^2}\\
&&\;\;\;\;\;\;\;\Big[\sum_{i,j=1}^l\bar{\alpha}_i\alpha_j\langle (Tf_i)^{k+1},(Tf_j)^{k+1}\rangle\langle B^{+(n-k)}_{Tf_i}\Phi,B^{+(n-k)}_{Tf_j}\Phi\rangle\Big].
\end{eqnarray*}
Put
$$M_k=\big(\langle f_i^{k+1},f_j^{k+1}\rangle\big)_{i,j},\;M_{k,T}=\big(\langle (Tf_i)^{k+1},(Tf_j)^{k+1}\rangle\big)_{i,j}.$$
This gives
\begin{eqnarray*}
\langle v,M_{k,T}v\rangle&=&\sum_{i,j=1}^l\bar{\alpha}_i\alpha_j\langle (Tf_i)^{k+1},(Tf_j)^{k+1}\rangle\\
&=&\|\alpha_1(Tf_1)^{k+1}+\dots+\alpha_l(Tf_l)^{k+1}\|^2_2\\
&=&\|\varphi^{k+1}T_1\big(\alpha_1f_1^{k+1}+\dots+\alpha_lf^{k+1}_l\big)\|^2_2\\
&\leq& \|\varphi\|^k_\infty\|T(\alpha_1f_1^{k+1}+\dots+\alpha_lf_l^{k+1})\|^2_2\\
&\leq&\|\alpha_1f_1^{k+1}+\dots+\alpha_lf_l^{k+1}\|^2_2=\langle v,M_kv\rangle.
\end{eqnarray*}
This proves that 
\begin{equation}\label{MKT}
0\leq M_{k,T}\leq M_k.
\end{equation}
Note that by induction assumption
\begin{equation}\label{An-k}
0\leq A_{n-k,T}\leq A_{n-k},
\end{equation}
for all $k=0,\dots,n$. Therefore, Lemma \ref{Kernel} and identies (\ref{MKT}), (\ref{An-k}) implies that
$$A_{n+1,T}\leq A_{n+1}.$$
Hence, we have proved that
$$\langle v,A_{n,T}v\rangle=\sum_{i,j=1}^l\bar{\alpha}_i\alpha_j\langle B^{+n}_{Tf_i}\Phi,B^{+n}_{Tf_j}\Phi\rangle\leq \langle v,A_{n}v\rangle=\sum_{i,j=1}^l\bar{\alpha}_i\alpha_j\langle B^{+n}_{f_i}\Phi,B^{+n}_{f_j}\Phi\rangle,$$
for all $n$. After using (\ref{chah}), it is easy to conclude that $\Gamma_2(T)$ is a contraction. 
\end{proof}

{\bf Remark}\\
The contractions considered in Proposition \ref{spec-contr} are very special,
however they are sufficient to prove the existence of the quadratic free Hamiltonian and the quadratic Ornstein--Uhlenbeck semigroup. In fact taking
$$
T=T=e^{z1_\mathcal{A}}
$$
with $Re(z)\leq0$ where $\mathcal{A}=L^2(\Bbb R^d)\cap L^\infty(\Bbb R^d)$,
Proposition \ref{spec-contr} implies that 
$\Gamma_2(e^{z1_\mathcal{A}})$ is a holomorphic semigroup which, by the 
remark done at the beginning of section (\ref{Quadr-sec-quant}), 
leaves the vacuum vector $\Phi$ invariant. 
In particular, for $z=it$, $t\in\mathbb{R}$, the generator $H_0$ 
of the strongly continuous $1$--parameter unitary group
$$
\Gamma_2(e^{it1_\mathcal{A}})=e^{it{H}_0}
$$ 
is the quadratic analogue of the free Hamiltonian.
By analytic continuation one has
$$
\Gamma_2(e^{z1_\mathcal{A}})=e^{z{H}_0}
$$
Moreover Lemma \ref{fir} shows that its action on the $n$--particle space
is the same as the action of the number operator in the usual Fock space, i.e.
it is reduced to multiplication by
$$
e^{zn}
$$
Thus $H_0$ is the positive self--adjoint operator characterized by the property that,
for any $n\in\mathbb{N}$, the $n$--particle space is the eigenspace of $H_0$ 
corresponding to the the eigenvalue $n$. \\
By considering the action of the {\it number operators} $N_f$,
defined at the beginning of section (\ref{Intr}), one easily verifies that
the definition of $N_f$ can be extended to the case in which $f$ is a multiple of the
identity function $1$, so that $N_1$ is well defined. With this notation one has 
the identity
$$
H_0 = {1\over 2} N_1
$$
Using the functional realization of the quadratic Fock space given by 
Theorem \ref{quad-IFS} it is clear that the contraction semigroup
$$
\Gamma_2(e^{-t1_\mathcal{A}})=e^{-t{H}_0}
$$
is positivity preserving and its explicit form gives that
 $$
\Gamma_2(e^{-t1_\mathcal{A}})1=e^{-t}1\leq 1
$$
(here we are extending in the obvious way the action of 
$\Gamma_2(e^{-t1_\mathcal{A}})$ to the multiples of the identity function
which is not in $L^2(\mathbb{R}^d)\cap L^\infty(\mathbb{R}^d))$. This means that the semigroup $e^{-t{H}_0}$ is sub--Markovian.
The above discussion shows that $e^{-t{H}_0}$ is a natural candidate 
for the role of quadratic analogue of the Ornstein--Uhlenbeck semigroup. 
A more detailed analysis of this semigroup
and of its properties will be discussed elsewhere.

\subsection{ A counterexample} 

In this subsection, we discuss the behavior of contractions under
quadratic second quantization.
\begin{lemma}\label{amna}
Let $T$ be a linear operator on $L^2(\Bbb R^d)\cap L^\infty(\Bbb R^d)$. If
 $T$ is a contraction on $L^2(\Bbb R^d)$ and on $L^\infty(\Bbb R^d)$, then for any quadratic exponential vector $\Psi(f)$ one has
\begin{eqnarray}\label{contr-QEV}
\|\Gamma_2(T)\Psi(f)\|\leq \|\Psi(f)\|
\end{eqnarray}
\end{lemma}
\begin{proof} 
Recall that
\begin{eqnarray}\label{norm-QEV}   
\|\Gamma_2(T)\Psi(f)\|^2=\|\Psi(Tf)\|^2=\sum_{n\geq0}\frac{\|B^{+n}_{Tf}\Phi\|^2}{(n!)^2}
\end{eqnarray}
and that, because of Lemma \ref{fir}:
\begin{eqnarray*}
\|B^{+n}_{Tf}\Phi\|^2=
\sum_{i_1+2i_2+\dots+ki_k=n}\frac{2^{2n-1}(n!)^2c^{i_1+\dots+i_k}}
{i_1!\dots i_k!2^{i_2}\dots k^{i_k}}\|Tf\|^{i_1}_2\|(Tf)^2\|^{i_2}_2\dots\|(Tf)^k\|^{i_k}_2
\end{eqnarray*}
If $T$ is a contraction on $L^2(\Bbb R^d)$ and on $L^\infty(\Bbb R^d)$  then
by the Riesz--Thorin Theorem, for all $p\geq2$, $T$ is also a contraction 
from $L^p(\Bbb R^d)$ into itself. Therefore, for any $p \geq 1$ and $i \in \mathbb N$:
$$
\|(Tf)^p\|^i_2=\left[\left(\int |Tf|^{2p}\right)^{1/2p}\right]^{pi}=\|Tf\|^{pi}_{2p} \leq 1
$$
for all $j=1,\dots,k$. This proves that for any $n\in\mathbb{N}$ 
$$
\|B^{+n}_{Tf}\Phi\|^2\leq\|B^{+n}_{f}\Phi\|^2
$$
and, in view of (\ref{norm-QEV}), this implies (\ref{contr-QEV}).
\end{proof}

From Lemma (\ref{amna}) it follows that the fact that $T$ is a contraction 
for $L^2(\mathbb{R}^d)$ and for $L^\infty(\mathbb{R}^d)$ is a necessary 
condition for $\Gamma_2(T)$ to be a contraction. 
The following counterexample shows that this condition is not sufficient.\\ \\
Define the linear operator 
$T:L^2(\mathbb{R})\cap L^\infty(\mathbb{R})\rightarrow 
L^2(\mathbb{R})\cap L^\infty(\mathbb{R})$ by
$$
Tf=\Big(\int_0^1f(t)dt\Big)\;\chi_{[0,1]}
$$
It is easy to verify that $T$ is a contraction in both $L^2$ and $L^\infty$. 
Therefore, from Lemma \ref{amna}, one has
$$
\|\Gamma_2(T)\Psi(f)\|\leq \|\Psi(f)\|
$$
In the following we will show that some linear combinations of
quadratic exponential vectors violate the inequality
$$
\|\Gamma_2(T)\Big(\sum_i\alpha_i\Psi(f_i)\Big)\|\leq\|\sum_i\alpha_i\Psi(f_i)\|
$$
In fact taking
$$
f_1:=\lambda\chi_{[0,\frac{1}{2}]}\quad ;\quad 
f_2:=\lambda\chi_{[0,1]}\quad ;\quad \lambda\in\mathbb{R}\quad ;\quad |\lambda|<\frac{1}{2}
$$
one has
$$
Tf_1=\frac{\lambda}{2}\chi_{[0,1]}\qquad ;\qquad Tf_2=\lambda\chi_{[0,1]}
$$
and  
\begin{eqnarray*}
\|\Gamma_2(T)\Big(\alpha_1\Psi(f_1)+\alpha_2\Psi(f_2)\Big)\|^2=\langle 
\left(
\begin{array}{lcc}
\alpha_1\\
\alpha_2
\end{array}
\right),B\left(
\begin{array}{lcc}
\alpha_1\\
\alpha_2
\end{array}
\right)\rangle,\\
\|\alpha_1\Psi(f_1)+\alpha_2\Psi(f_2)\|^2=\langle 
\left(
\begin{array}{lcc}
\alpha_1\\
\alpha_2
\end{array}
\right),A\left(
\begin{array}{lcc}
\alpha_1\\
\alpha_2
\end{array}
\right)\rangle
\end{eqnarray*}
where the matrices $A,B$ are defined by:
$$
A:=(\langle\Psi(f_i),\Psi(f_j)\rangle)_{1\leq i,j\leq2}
\qquad ;\qquad B:=(\langle\Psi(Tf_i),\Psi(Tf_j)\rangle)_{1\leq i,j\leq2}
$$
The contraction condition
$$
\|\Gamma_2(T)\Big(\alpha_1\Psi(f_1)+\alpha_2\Psi(f_2)\Big)\|^2\leq \|\alpha_1\Psi(f_1)+\alpha_2\Psi(f_2)\|^2
$$
is equivalent to say that $B\leq A$. In the following we prove that this 
inequality is not true. In fact recalling (\ref{Form}), i.e.
$$
\langle \Psi(f),\Psi(g)\rangle
=e^{-\frac{c}{2}\int_\mathbb{R}\ln(1-4\bar{f}(x)g(x))dx}
$$
one finds
$$
A=\left(
\begin{array}{lcc}
(\frac{1}{1-4\lambda^2})^{\frac{c}{4}} & (\frac{1}{1-4\lambda^2})^{\frac{c}{4}}\\
(\frac{1}{1-4\lambda^2})^{\frac{c}{4}} & (\frac{1}{1-4\lambda^2})^{\frac{c}{2}}
\end{array}
\right),\;
B=\left(
\begin{array}{lcc}
(\frac{1}{1-\lambda^2})^{\frac{c}{2}} & (\frac{1}{1-2\lambda})^{\frac{c}{2}}\\
(\frac{1}{1-2\lambda^2})^{\frac{c}{2}} & (\frac{1}{1-4\lambda^2})^{\frac{c}{2}}
\end{array}
\right)
$$
and a simple calculation proves that $det(A-B)\leq0$. \\
\bigskip\\

{\bf\large Acknowledgments}\bigskip

Ameur Dhahri gratefully acknowledges stimulating discussions with Uwe Franz and Eric Ricard.

\end{document}